\documentclass[12pt]{article}
\usepackage{graphicx,color,amsmath,amssymb}
\newtheorem{theorem}{Theorem}\newtheorem{lemma}{Lemma}

\newtheorem{claim}{Claim}
\newtheorem{proposition}{Proposition}

\def\beq{\begin{equation}}\def\eeq{\end{equation}}
\def\beqn{\begin{eqnarray}}\def\eeqn{\end{eqnarray}}
\def\pont{\hspace{-6pt}{\bf.\ }}
\def\eps{\varepsilon}
\def\qed{\ifhmode\unskip\nobreak\fi\quad\ifmmode\Box\else$\Box$\fi}
\def\ex{{\rm ex}}
\title{Tur\'an number of special four cycles in triple systems}

\author{
Zolt\'an F\"uredi, Andr\'as Gy\'arf\'as, Attila Sali\\
\small Alfr\'ed R\'enyi Institute of Mathematics\\[-0.8ex]
\small Hungarian Academy of Sciences\\[-0.8ex]}

\begin{document}
\maketitle

\begin{abstract} A {\em special four-cycle } $F$ in a triple system consists of four triples {\em inducing } a $C_4$. This means that $F$ has four special vertices $v_1,v_2,v_3,v_4$ and four triples in the form $w_iv_iv_{i+1}$ (indices are understood $\pmod 4$) where the $w_j$s are not necessarily distinct but disjoint from $\{v_1,v_2,v_3,v_4\}$. There are seven non-isomorphic special four-cycles, their family is denoted by $\cal{F}$. Our main result implies that the Tur\'an number $\ex(n,{\cal{F}})=\Theta(n^{3/2})$. In fact, we prove more, $\ex(n,\{F_1,F_2,F_3\})=\Theta(n^{3/2})$, where the $F_i$-s are specific members of $\cal{F}$. This extends previous bounds for the Tur\'an number of triple systems containing no Berge four cycles.

We also study $\ex(n,{\cal{A}})$ for all ${\cal{A}}\subseteq {\cal{F}}$. For 16 choices of $\cal{A}$ we show that $ex(n,{\cal{A}})=\Theta(n^{3/2})$, for 92 choices of $\cal{A}$ we find that $\ex(n,{\cal{A}})=\Theta(n^2)$ and the other 18 cases remain unsolved.

\end{abstract}

\section{Introduction}

A {\em triple system } $H=(V,E)$ has vertex set $V$ and $E$ consists of some triples of $V$ (repeated triples are excluded). For any fixed family  ${\cal{H}}$ of triple systems, the Tur\'an number $\ex(n,{\cal{H}})$ is the maximum number of triples in a triple system of $n$ vertices that is $\cal{H}$-free, i.e., does not contain any member of ${\cal{H}}$ as a subsystem.

Our interest here is the family ${\cal{F}}$ of {\em special four cycles}: they have four distinct base vertices $v_1,v_2,v_3,v_4$ and four triples $w_iv_iv_{i+1}$ (indices are understood $\pmod 4$) where the $w_j$s are not necessarily distinct but $w_i\ne v_j$ for any pair of indices $1\le i,j\le 4$.

\begin{figure}
\begin{center}
\includegraphics[scale=0.5]{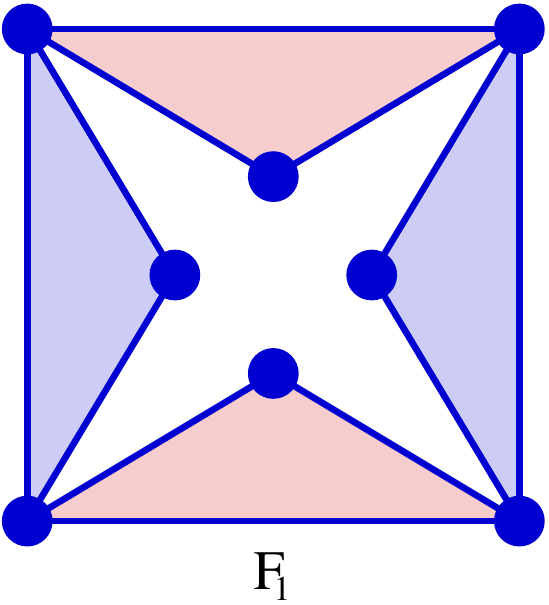}\hfill\includegraphics[scale=0.5]{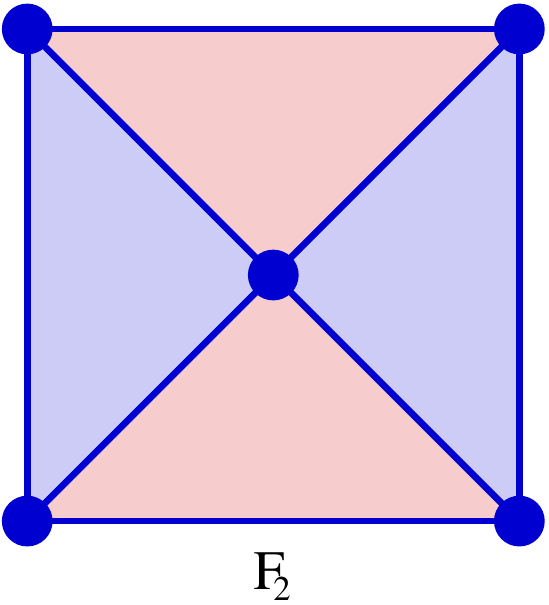}\hfill\includegraphics[scale=0.5]{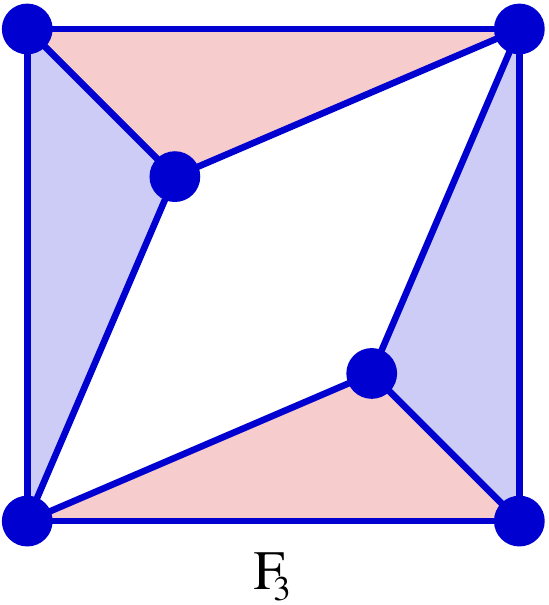}
\end{center}
\caption{The family of special four cycles $F_1,F_2,F_3$ of Theorem~\ref{main}}\label{f1f3}
\end{figure}

\begin{figure}
\begin{center}
\includegraphics[scale=0.5]{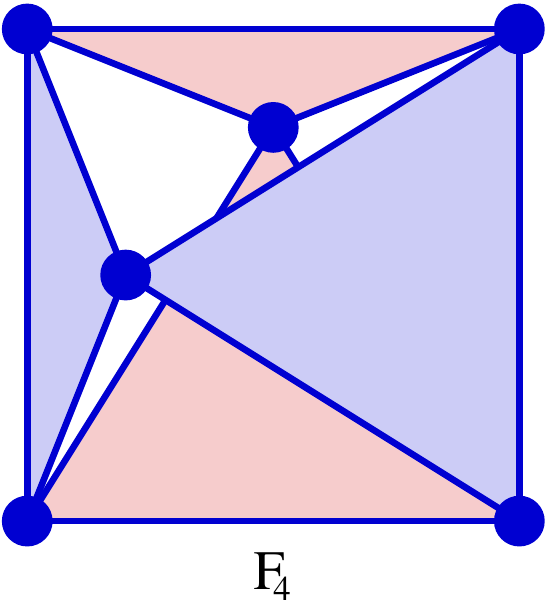}\hfill\includegraphics[scale=0.5]{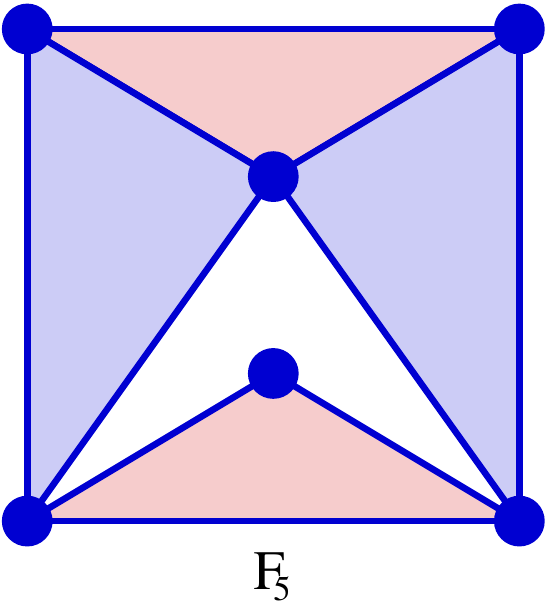}\hfill\includegraphics[scale=0.5]{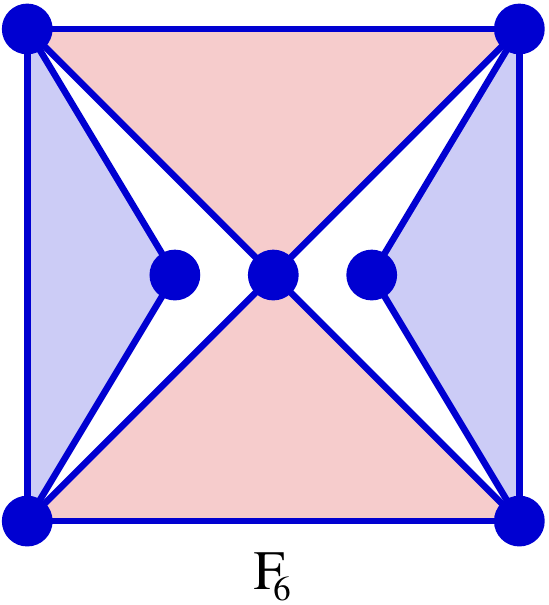}\hfill\includegraphics[scale=0.5]{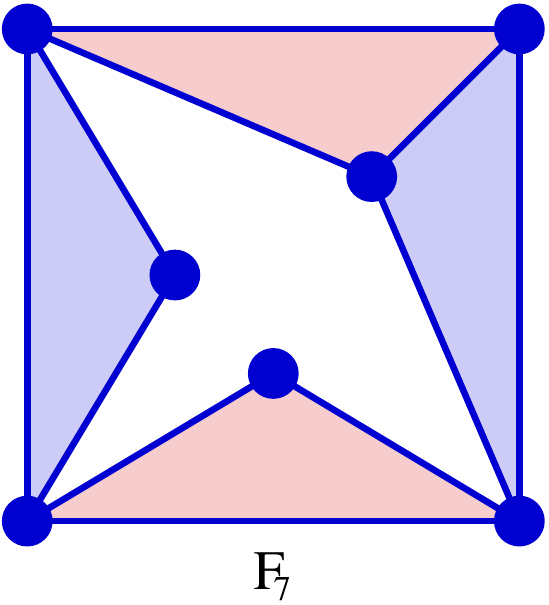}
\end{center}
\caption{The family of the other four special four cycles $F_4,\dots, F_7$}\label{f4f7}
\end{figure}

There are seven non-isomorphic special four cycles. The linear (loose) four cycle $F_1$ is obtained when all $w_j$-s are different and in $F_2$ all $w_j$s coincide. When two pairs coincide we get either $F_3$ ($w_1=w_2,w_3=w_4$) or $F_4$ ($w_1=w_3,w_2=w_4$). The $F_4$ is the Pasch configuration. We define $F_5$ with $w_1=w_2=w_3$ (but $w_4$ is different). In $F_6$ we have $w_1=w_3$ (and $w_2,w_4$ are different from $w_1$ and from each other).  When only $w_1,w_2$ coincide we get $F_7$. Set ${\cal{F}}=\{F_1,\dots,F_7\}$. For the convenience of the reader, the special four cycles are shown on Figure~\ref{f1f3} and Figure~\ref{f4f7}.

Tur\'an numbers of various members of ${\cal{F}}$ have been investigated before. F\"uredi \cite{FU1} proved that $\ex(n,F_3)\le {7\over 2}{n\choose 2}$. Mubayi \cite{MU} showed that $\ex(n,F_2)=\Theta(n^{5/2})$. R\"odl and Phelps \cite{ROPH} gave the bounds   $c_1n^{5/2}\le \ex(n,F_4)\le c_2n^{11/4}$. In fact, the upper bound is Erd\H os' upper bound \cite{E} for $\ex(n,K^3_{2,2,2})$. The lower bound comes from a balanced $3$-partite triple system where every vertex of the third partite class form a triple with the edges of a bipartite $C_4$-free graph between the first two partite classes.

We prove that $\ex(n,{\cal{F}})=\Theta(n^{3/2})$, thus has the same order of magnitude as $\ex(n,C_4)$ for graphs. In fact, it is enough to exclude three of the special four cycles.

\begin{theorem}\pont \label{main} $\ex(n,\{F_1,F_2,F_3\})=\Theta(n^{3/2})$.
\end{theorem}

The family ${\cal{F}}$ of special four cycles is a subfamily of a wider class, the class of {\em Berge four cycles}, where the vertices $w_i$ can be selected from the base vertices as well, requiring only that the four triples $w_iv_iv_{i+1}$ are different. Theorem \ref{main} extends previous similar upper bounds (F\"uredi and \"Ozkahya \cite{FO}, Gerbner, Methuku, Vizer \cite{GMV}) where the family of Berge four cycles were forbidden.

The appearance of the set $\{F_1,F_2,F_3\}$ is not accidental. If any of $F_1,F_2,F_3$ is missing from ${\cal{A}}\subset {\cal{F}}$ then $\ex(n,{\cal{A}})$ is essentially larger than $n^{3/2}$.

\begin{itemize}\label{constr}
\item (C1) Ruzsa and Szemer\'edi \cite{RUSZ} constructed triple systems on $n$ vertices that do not carry three triples on six vertices and have more than $n^{2-\eps}$ triples for any fixed $\eps$. This provides an example which contains only $F_1$ from ${\cal{F}}$,

\item (C2) The ${n-1\choose 2}$ triples containing a fixed vertex from $n$ vertices contains only $F_2$ from ${\cal{F}}$,

\item (C3) Partition $n$ vertices evenly into three parts, take a pairing between two equal parts and extend each pair with all vertices of the third class to a triple. This gives a triple system with approximately $n^2/9$ triples and contains only $F_3$ from ${\cal{F}}$.

\end{itemize}

In Section 3 we discuss $\ex(n,{\cal{A}})$ for all ${\cal{A}}\subseteq {\cal{F}}$. It turns out that in 92 cases $\ex(n,{\cal{A}})=\Theta(n^2)$ and 18 cases remain unsolved.
\section{Proof of Theorem \ref{main}}

Assume $H$ is a triple system with $n$ vertices containing no subsystem from the set $F_1,F_2,F_3$. Applying the standard approach (based on \cite{EK}), we may assume that $H$ is $3$-partite with vertex partition $[A_1,A_2,A_3]$  where $|A_i|\in \{\lfloor n/3 \rfloor, \lceil n/3   \rceil\}$ and contains at least $2/9$ of the triples of the original triple system.

The triples of $H$ define a bipartite graph $B=[A_1,A_2]$ as follows. If $(a_1,a_2,a_3)$ is a triple of $H$ with $a_i\in A_i$ then $a_1a_2$ is considered as an edge of $B$. Define the label $L(a_1,a_2)$ of $a_1a_2\in E(B)$ as the set $\{z\in A_3: a_1a_2z\in E(H)\}$. Then

\beq \label{size} |E(H)|=\sum_{a_1a_2\in E(B)} |L(a_1,a_2)|. \eeq


\begin{lemma}\pont \label{simple} The bipartite graph $B$ has at most $O(n^{3/2})$ edges.
\end{lemma}

\noindent {\bf Proof of Lemma \ref{simple}.} We denote by $N(x,y)$ the set of common neighbors (in $B$) of $x,y\in A_2$ in $A_1$. Similarly, let $N(u,v)$ be the set common neighbors of $u,v\in A_1$ in $A_2$.

For distinct vertices $x,y\in A_2$, define the digraph $D=D(x,y)$ with vertex set $A_3$.  For every $u\in A_1$ such that $u\in N(x,y)$ and $a_i\in L(u,x), a_j\in L(u,y)$, a directed edge $a_ia_j$ is defined in $D(x,y)$.  We claim that $D(x,y)$ is a very special digraph.
\begin{claim}\pont\label{spec}
\begin{itemize}
\item (1.1) There are no multiple loops or parallel directed edges in $D(x,y)$,
\item (1.2) There is at most one loop in $D(x,y)$,
\item (1.3) Two non-loop edges of $D(x,y)$  either intersect or $|N(x,y)|\le 4$.
\end{itemize}
\end{claim}
\noindent {\bf Proof.} A multiple loop $a_ia_i$ in $D(x,y)$ would give a $C_4$ in $B$ with all edges containing $a_i$ in their labels, this corresponds to an $F_2$ in $H$ -- a contradiction. A multiple edge $a_ia_j$ would give a $C_4=(x,u_1,y,u_2)$ in $B$ where $u_1,u_2\in N(x,y), u_1\ne u_2$ such that the  consecutive edges of $C_4$ contain $a_i,a_j,a_j,a_i\in A_3$ in their labels. This would give an $F_3$ in $H$ -- a contradiction again, proving (1.1).

Two distinct loops $a_ia_i,a_ja_j$ in $D(x,y)$ can appear in two ways. Either we have a $C_4=(x,u_1,y,u_2)$ in $B$ where $u_1,u_2\in N(x,y), u_1\ne u_2$ such that the consecutive edges of $C_4$ contain $a_i,a_i,a_j,a_j\in A_3$ in their labels, this would give an $F_3$ in $H$, a contradiction. Otherwise $u=u_1=u_2$ and we have two multiedges $xu,yu$ both containing  $a_i,a_j$ in their labels, this gives an $F_2$ in $H$ with $u$ in its center, a contradiction again, proving (1.2).

Suppose that there exists two non-intersecting non-loop edges $a_ia_j,a_ka_l$ in $D(x,y)$. If these edges are defined by $u_1,u_2\in N(x,y), u_1\ne u_2$, we have a $C_4=(x,u_1,y,u_2)$ in $B$ with four distinct elements in their labels, giving an $F_1$ in $H$, a contradiction. Thus we may assume that $u_1=u_2=u$ and we have $xu,yu$ in $B$ with $a_i,a_k$ and with $a_j,a_l$ in their labels. Set
$$M=\{v\in N(x,y): v\ne u, |L(v,x)\cup L(v,y)|\ge 2\}.$$ We claim that $|M|\le 2$. Indeed, consider $v\in M$, there is $a_s,a_t\in A_3$ such that $a_s\ne a_t$ and $xv,yv$ have labels containing $a_s,a_t$, respectively.
Observe that either $\{s,t\}=\{i,k\}$ or $\{s,t\}=\{j,l\}$ otherwise there is a $C_4=(x,u,y,v)$ with four distinct labels, giving an $F_1$ in $H$, a contradiction. This implies that $|M|\le 4$. However, it cannot happen that for two distinct vertices $v,v'\in M$ the coincidence of the index pairs are $\{i,k\}$ and $\{j.l\}$, respectively, because it would result again in a $C_4=(x,v,y,v')$ with four distinct labels, a contradiction as above.
Thus $|M|\le 2$ (equality is possible with edge pairs $a_ia_k,a_ka_i$ or $a_ja_l,a_la_j$), proving the claim.

Observing that every vertex of $N(x,y)\setminus (\{u\}\cup M)$ defines a loop in $D(x,y)$, (1.1) and  (1.2) implies that $|N(x,y)|\le 4$, proving  (1.3) and Claim \ref{spec}. \qed

A {\em cherry } on $x\in A_2$ is defined as an incident edge pair, $ux,vx\in E(B)$ such that $u,v\in A_1,u\ne v$ and $L(u,x)\cap L(v,x)\ne \emptyset$.  Let $C(x,y)$ be the number of cherries in the subgraph of $B$ induced on $\{x,y\}\cup N(x,y)$.
We claim
\beq \label{cxy} C(x,y)\ge \sum_{a\in V(D(x,y))} d^+(a)+d^-(a)-2,  \eeq
because there are atleast
$d^+(a)-1$ cherries on $x$ with $L(u,x)\cap L(v,x)=\{a\}$ and $d^-(a)-1$ cherries on $y$ with $L(u,y)\cap L(v,y)=\{a\}$.

\begin{claim}\pont\label{lowb} For any two distinct vertices $x,y\in A_2$, $C(x,y)\ge |N(x,y)|-4.$
\end{claim}

\noindent {\bf Proof.} It is certainly true  for $|N(x,y)|\le 4$. Otherwise, using (1.3) from Claim \ref{spec}, we have pairwise intersecting edges in $D(x,y)$.

\noindent {\bf Case 1.} The edges of $D(x,y)$ form a triangle (edges oriented two ways are allowed) plus at most one loop. Therefore $D(x,y)$ has at most seven edges thus $5\le |N(x,y)|\le 7$. 
by (1.1) of Claim~\ref{spec}, $d^+(a)\le 2$ for any vertex of the triangle. There are at least $|N(x,y)|-1$ edges on the triangle, so there exists at least $|N(x,y)|-4$ vertices $a$ with  $d^+(a)\ge 2$ resulting in at least  $|N(x,y)|-4$ cherries on $x$.

\noindent {\bf Case 2.} All edges of $D(x,y)$ (apart from a possible loop)  contain $a\in A_3$.  For every $u\in N(x,y)$ (apart from one possible vertex which defines a loop) either $ux$ or $uy$ has label $a$. Thus $\sum_{a\in V(D(x,y))} d^+(a)+d^-(a)\ge |N(x,y)|-1$  so (\ref{cxy}) results in at least $|N(x,y)|-3$  cherries on $x$ or on $y$, completing the proof of Claim \ref{lowb}. \qed

\begin{claim}\pont\label{ub} $\sum_{x,y\in A_2} C(x,y)\le {|A_1|\choose 2}$.
\end{claim}

\noindent {\bf Proof.} Every cherry counted on the left hand side is on some pair of $A_1$. At most one cherry can be on any $(u,v)\in A_1$, otherwise (by (1.1) in Claim \ref{spec}) we have one of $F_2,F_3$.  \qed

\smallskip

Applying Claims \ref{lowb}, \ref{ub}  we get
$$\sum_{x,y\in A_2}\left(|N(x,y)| -4\right) \le  \sum_{x,y\in A_2} C(x,y)\le {|A_1|\choose 2},$$
thus
$$\sum_{x,y\in A_2} |N(x,y)|\le 4{|A_2|\choose 2}+{|A_1|\choose 2}\le O(n^{2}).$$

By convexity we get

$$|A_1|{{|E(B)|\over |A_1|}\choose 2} \le \sum_{u\in A_1} {d(u)\choose 2}=\sum_{x,y\in A_2} |N(x,y)|\le O(n^2),$$
therefore $|E(B)|=O(n^{3/2})$, proving Lemma \ref{simple}.   \qed

\bigskip

To finish the proof of Theorem \ref{main}, we need to show that the presence of labels does not affect strongly the edge count of Lemma \ref{simple}. Let $B^*$ denote the subgraph of $B$ with the edges of at least three-element labels.

\begin{proposition}\pont\label{multi} If $H$ is $\{F_1,F_2,F_3\}$-free then $B^*$ is $C_4$-free.
\end{proposition}

\noindent {\bf Proof.}  Assume $C=(x,u,y,v,x)$ is a four-cycle in $B^*$. From the definition of $B^*$ there are  three distinct elements, say $a,b,c$ from the labels of three edges of $C$. The only way to avoid $F_1$ is that the fourth edge has label $\{a,b,c\}$. However, the same argument forces that all labels on $C$ are equal to $\{a,b,c\}$ giving (many) $F_3$'s. \qed


\smallskip

We can consider $B^*$ as a bipartite multigraph obtained as the union of $|A_3|$ simple bipartite graphs as follows. Set
$$E(z)=\{(u,x): u\in A_1, x\in A_2,  (u,x,z)\in E(H)\text{ and } |L(u,x)|\ge 3\},$$
then $E(B^*)=\cup_{z\in A_3} E(z)$.

\begin{proposition}\pont\label{nop5} For every $z\in A_3$ there is no path in $B^*$ with four edges such that its first and last edge is in $E(z)$.
\end{proposition}

\noindent {\bf Proof.}   Suppose that edges $e_1,e_2,e_3,e_4$ form such a path for some $z\in A_3$. Since each edge of $B^*$ has multiplicity at least three, we can replace $e_2$ by $f_2$ and $e_3$ by $f_3$ so that $f_2\in E(z_1), f_3\in E(z_2)$ and $z_1,z_2$ are distinct and both different from $z$. Then the four triples of $H$, $$e_1\cup \{z\},f_2\cup \{z_1\},f_3\cup \{z_2\},e_4\cup \{z\}$$
form an $F_1$, contradiction.
\qed

\smallskip

For any vertex $x\in A_2$ let $L(x)$ denote the subset of $A_3$ that appears in some of the labels on edges of $B^*$ incident to $x$.

\begin{proposition}\pont\label{noreppair} For distinct vertices $x_1,x_2,x_3,x_4 \in A_2$,
$$|L(x_1)\cap L(x_2)\cap L(x_3)\cap L(x_4)|\le 1.$$
\end{proposition}

\noindent {\bf Proof.}  Suppose on the contrary that we have $z_1,z_2\in A_3$ such that for $i=1,2,3,4$, $e_i=\{z_1,x_i,u_{2i-1}\}, f_i=\{z_2,x_i,u_{2i}\}$ are all triples of $H$.

An $F_1$ is formed by the triples $e_i,f_i,e_j,f_j$ if there is a pair $i,j$ such that \\$u_{2i-1},u_{2i},u_{2j-1},u_{2j}$ are all different. Thus, we may assume that for any pair $1\le i<j\le 4$ there is an equality between elements  $u_{2i-1},u_{2i},u_{2j-1},u_{2j}$.

Let us call an equality $u_{2i-1}=u_{2i}$ \emph{horizontal}, an equality $u_{2i}=u_{2j}$ or $u_{2i-1}=u_{2j-1}$ (for $i\ne j$) \emph{vertical}, finally an equality  $u_{2i-1}=u_{2j}$ (for $i\ne j$) \emph{diagonal}. The terms to distinguish  equalities refer to an arrangement of the vertices $u_i$ into a $4\times 2$ matrix with $u_{2i-1}, u_{2i}$ in row $i$. Observe the following facts.
\begin{enumerate}
\item \label{label:2horiz}$F_3$ or $F_2$ is formed by the triples $e_i,f_i,e_j,f_j$ if the pair $i\ne j$ have both horizontal equalites holding. Thus, at most one horizontal equality may hold.
\item \label{label:2vert} If there is pair  $i\ne j$ such that both vertical equalities hold, then a $C_4$ can be found in $B^*$ contradicting to Proposition~\ref{multi}. Similarly,
  \item \label{label:2diag} if there is pair  $i\ne j$ such that both diagonal equalities hold, we get a contradiction with Proposition~\ref{multi}.

\item \label{label:1vert} We get a four edge path contradicting to Proposition \ref{nop5} if there are is a pair $i,j$ such that exactly one vertical equality holds, that is  $u_{2i}=u_{2j}$ and $u_{2i-1},u_{2j-1}$ are different and different from $u_{2i}$ as well. (Symmetrically, if there are $x_i,x_j$ such that $u_{2i-1}=u_{2j-1}$ and $u_{2i},u_{2j}$ are different and different from $u_{2i-1}$ as well.)
\end{enumerate}

Facts \ref{label:2horiz}--\ref{label:1vert} imply that there exists a triple of indices, $i,j,k$ such that we have exactly one diagonal equality on each pair of them. These are either in the form $u_{2i-1}=u_{2k},u_{2i}=u_{2j-1},u_{2j}=u_{2k-1}$
defining a six-cycle in $B$ on the vertices $x_i,x_j,x_k,u_{2i},u_{2j},u_{2k}$, giving  (three) $F_1$, for example $e_i,f_k,f_j,e_j$, or in the form $u_{2i-1}=u_{2k},u_{2i}=u_{2j-1},u_{2j-1}=u_{2k}$ that implies horizontal equality $u_{2i-1}=u_{2i}$, a contradiction.
This proves Proposition \ref{noreppair}.\qed

\smallskip

By Propositions \ref{multi}, \ref{nop5}  the simple bipartite graph $B(z)$ with edge set $E(z)$ has no cycles or paths with four edges. Therefore each component of $B(z)$ is a double star. Thus each $B(z)$ can be written as the union of two graphs, $S(z),T(z)$ where each vertex of $S(z)\cap A_2$ and each vertex of $T(z)\cap A_1$  has degree one in $B(z)$.  Set $$S=\cup_{z\in A_3} S(z), T=\cup_{z\in A_3} T(z).$$
By the definition of $S$, for every vertex $x\in A_2$, we have $|L(x)|=d_S(x)$ where $d_S(x)$ is the degree of vertex $x$ in the (multi) graph
$S$. By Proposition \ref{noreppair} $$\sum_{x\in A_2} {|L(x)|\choose 2}\le 3{|A_3|\choose 2}$$ therefore $$\sum_{x\in A_2} {d_S(x)\choose 2}\le 3{|A_3|\choose 2}.$$

Applying the same argument symmetrically for vertices of $A_1$ and for the graph $T$, we get
$$\sum_{u\in A_1} {d_T(u)\choose 2}\le 3{|A_3|\choose 2}.$$
By the convexity argument, $|E(B^*)|=|E(S)|+|E(T)|=O(n^{3/2})$. By Lemma \ref{simple}, we also have $|E(B)|=O(n^{3/2})$. Thus by (\ref{size}) and the definition of $B^*$,
$$|E(H)|=\sum_{a_1a_2\in E(B)} |L(a_1,a_2)|\le 2|E(B)|+|E(B^*)|=O(n^{3/2}),$$
concluding the proof of Theorem \ref{main}.

\section{Concluding remarks}

Theorem \ref{main} determines the order of magnitude ($\Theta(n^{3/2})$) for the 16 subsets of $\cal{F}$ containing $F_1,F_2,F_3$ and we pointed out that for all other choices ${\cal{A}}\subset {\cal{F}}$, $\ex(n,{\cal{A}})$ must be essentially larger. In this section we summarize what we know about these cases. There is a trivial case, when $\cal{A}$ is empty and $\ex(n,{\cal{A}})={n\choose 3}$. Furthermore, as mentioned before, $\ex(n,F_2)=\Theta(n^{5/2})$ was proved by the first author (see in Mubayi \cite{MU}). Thus we have $2^7-2^4-2=110$ cases to consider. It turns out that in 92 cases the order of magnitude is $\Theta(n^2)$ (see Subsection \ref{quadratic}) and only the remaining 18 cases are left unsolved (see Subsection \ref{unsolved}).

A simple but useful lemma compares Tur\'an numbers of closely related triple systems. Assume $G$ is a triple system and $v,w\in V(G)$ is covered by $e\in E(G)$. The triple system obtained from $G$ by removing $e$ and adding the triple $v,w,x$ where $x\notin V(G)$ is called a {\em fold out} of $G$. For example $F_7$ is a fold out of $F_3$, $F_6$ is a fold out of $F_4$.

\begin{lemma}\pont\label{fold} (Fold out lemma.) If $G$ is a triple system and $G_1$ is a fold out of $G$ then $\ex(n,G_1)\le ex(n,G)+(|V(G)|-2){n\choose 2}$.
\end{lemma}

\noindent {\bf Proof.} Suppose that a triple system $H$ has $n$ vertices and has more than $\ex(n,G)+(|V(G)|-2){n\choose 2}$ triples. A triple of $H$ is called bad if it contains a pair of vertices that covered by at most  $|V(G)|-2$ triples of $H$, otherwise it is a good triple. Then $H$ has more than $\ex(n,G)$ good triples thus contains a copy of $G$ with all triples good. By definition, any pair of vertices in any triple of this copy of $G$ is in more than $|V(G)|-2$ triples of $H$ so some of them is suitable to define the required fold out $G_1$ of $G$.  \qed

\subsection{ When $\ex(n,{\cal{A}})=\Theta(n^2)$}\label{quadratic}

Here we collect all cases of ${\cal{A}}\subset {\cal{F}}$ when we can prove that $\ex(n,{\cal{A}})=\Theta(n^2)$.

\begin{proposition}\pont\label{nsq1} Assume that ${\cal{A}} \subset {\cal{F}}\setminus F_2$ and ${\cal{A}} \cap \{F_1,F_3,F_7\}\ne \emptyset$. Then $\ex(n,{\cal{A}})=\Theta(n^2)$.
\end{proposition}

\noindent {\bf Proof.} The first condition ensures that the members of ${\cal{A}}$ cannot be pierced by one vertex, thus
Construction~(C3) shows that $\ex(n,{\cal{A}})=\Omega(n^2)$. On the other hand, $F_7$ is a fold out of $F_3$ and $F_1$ is a fold out of $F_7$ thus by Lemma \ref{fold} (and by the second condition of the proposition)
$$\ex(n,F_1)\le \ex(n,F_7)+O(n^2)\le \ex(n,F_3)+O(n^2)\le {7\over 2}{n\choose 2}+O(n^2)$$ where the upper bound of $\ex(n,F_3)$ is F\"uredi's result \cite{FU1}. \qed

\begin{proposition}\pont \label{nsq2} Assume that ${\cal{A}} \subset {\cal{F}}\setminus F_3$ and ${\cal{A}} \cap \{F_1,F_7\}\ne \emptyset$. Then $\ex(n,{\cal{A}})=\Theta(n^2)$.
\end{proposition}

\noindent {\bf Proof.} To show that $\ex(n,F)=\Omega(n^2)$, consider the Construction (C3), it contains only $F_3$ from $\cal{F}$.  The upper bound follows by the argument of Proposition \ref{nsq1}. \qed

\begin{proposition}\pont \label{nsq3} Assume that $\{F_2,F_3\} \subset {\cal{A}} \subset \{F_2,F_3,F_4,F_5,F_7\}$. Then $\ex(n,{\cal{A}})=\Theta(n^2)$.
\end{proposition}

\noindent {\bf Proof.} To show that $\ex(n,{\cal{A}})=\Omega(n^2)$, consider
\begin{itemize}
\item (C4) Steiner triple systems without $F_4$ (the Pasch configuration), they do not contain any member of $\cal{A}$.
\end{itemize}
The upper bound follows from \cite{FU1} since $F_3\in {\cal{A}}$.  \qed

\begin{proposition}\pont \label{nsq4} Assume that $\{F_2,F_3,F_6\} \subset {\cal{A}} \subset \{F_2,F_3,F_5,F_6,F_7\}$. Then $\ex(n,{\cal{A}})=\Theta(n^2)$.
\end{proposition}

\noindent {\bf Proof.} To show that $\ex(n,{\cal{A}})=\Omega(n^2)$, consider
\begin{itemize}
\item (C5) Steiner triple systems without $F_6$ (projective Steiner triple systems), they do not contain any member of $\cal{A}$.
\end{itemize}
The upper bound follows again from \cite{FU1} since $F_3\in {\cal{A}}$. \qed

Note that Proposition \ref{nsq1} covers 56 cases, Proposition \ref{nsq2} adds 24 further cases, Propositions \ref{nsq3}, \ref{nsq4} add 8 plus 4 further cases. These 92 cases are the ones when  $\ex(n,{\cal{A}})=\Theta(n^2)$ follows from known results.

\subsection{Unsolved cases}\label{unsolved}

The 18 unsolved cases are grouped as follows.
\begin{itemize}
\item  1. $\ex(n,\{{\cal{A}}\cup F_6\})$ where  ${\cal{A}}\subseteq \{F_2,F_4,F_5\}$ ( 8 cases )

\item  2. $\ex(n,\{F_2,F_5\}), \ex(n,F_5)$

\item  3. $\ex(n,\{F_2,F_4,F_5\}), \ex(n,\{F_4,F_5\})$

\item  4. $\ex(n,\{F_2,F_4\}), \ex(n,F_4)$

\item 5. $\ex(n,{\cal{A}})$ where $\{F_2,F_3,F_4,F_6\}\subseteq {\cal{A}} \subseteq \{F_2,F_3,F_4,F_5,F_6,F_7\}$ (4 cases)
\end{itemize}

The upper bounds for the unknown cases can be compared by using Lemma \ref{fold}. For example, observing that $F_6$ is a fold out of $F_4$ and of $F_5$, moreover $F_5$ is a fold out of $F_2$, Lemma \ref{fold} implies

\begin{proposition}\pont\label{f6} Let $\cal{A}$ be any subset of $\{F_2,F_4,F_5\}$. Then $$\ex(n,F_6) \leq \ex(n,\{{\cal{A}}\cup F_6\})\le \ex(n,F_6)+7{n\choose 2}.$$
\end{proposition}

A lower bound $\Omega (n^2)$ for the first four groups of unknown cases can be obtained from Construction~(C3).
Lower bounds for the fifth group of unknown cases can be given by well studied functions introduced in \cite{BES}. Let $\ex(n,(6,3))$ be the maximum number of triples in a triple system that does not contain three triples inside any six vertices. Since all members of ${\cal{F}}$ except $F_1$ contain three triples inside six vertices an almost quadratic lower bound of Construction (C1) comes from~\cite{RUSZ} for the four unsolved cases in group 5. A quadratic upper bound is from \cite{FU1} since $F_3\in {\cal{A}}$. Thus we get

\begin{proposition}\pont\label{lb} If $\{F_2,F_3,F_4,F_6\}\subseteq {\cal{A}} \subseteq \{F_2,F_3,F_4,F_5,F_6,F_7\}$ then $$\ex(n,(6,3))\le \ex(n,{\cal{A}})=O(n^2)$$.
\end{proposition}
In fact, the lower bound of Proposition \ref{lb} can be changed to $\ex(n,(7,4))$ (the maximum number of triples in a triple system that does not contain four triples inside any seven vertices).

\begin{paragraph}{Acknowledgement}
  The authors are indebted to an unknown referee for pointing out some inaccuracies in the manuscript.

{\small  The research of the first author was supported in part by the Hungarian National Research, Development and Innovation Office NKFIH grant KH-130371 and NKFI–133819.}
\end{paragraph}


\begin{thebibliography}{99}
\bibitem{BES} W. G. Brown, P. Erd\H os, V. T. S\'os, Some extremal problems on $r$-graphs, in {\em New directions in the theory of graphs, Proc. 3rd Ann Arbor Conference on Graph Theory}, Academic Press, New York, 1973, 55--63.
\bibitem{E} P. Erd\H os, On extremal problems of graphs and generalized graphs, {\em Israel Journal of Mathematics} {\bf 2} (1964) 183--190.
\bibitem{EK} P. Erd\H os, D. J. Kleitman, On coloring graphs to maximize the proportion of of multicolored $k$-edges, {\em Journal of Combinatorial Theory } {\bf 5 } (1968) 164--169.
\bibitem{FU1} Z. F\"uredi, Hypergraphs in which all disjoint pairs have distinct unions, {\em Combinatorica} {\bf 4} (1984) 161--168.
\bibitem{FO} Z. F\"uredi, L. \"Ozkahya, On 3-uniform hypergraphs without a cycle of given length, {\em Discrete Applied Mathematics} {\bf 216} (2017), 582--588.
\bibitem{GMV} D. Gerbner, A. Methuku, M. Vizer, Asymptotics for the Tur\'an number of Berge-$K_{2,t}$, arXiv:1705.04134v2
\bibitem{MU} D. Mubayi, Some exact results and new asymptotics for hypergraph Tur\'an numbers, {\em Combinatorics, Probability and Computing} {\bf 11} (2002) 299--309.
\bibitem{ROPH}H. Leffmann, K. T. Phelps, V. R\"odl, Extremal problems for triple systems, {\em Journal of Combinatorial Designs} {\bf 1.} (1993) 379--394.
\bibitem{RUSZ} I. Z. Ruzsa, E. Szemer\'edi, Triple systems with no six points carrying three triangles, {\em in: Combinatorics, Vol. II. Coll. Math. Soc. J. Bolyai} {\bf 18 } North-Holland, 1978, 939--945.

\end{thebibliography}
\end{document}